\documentclass{amsart}

\usepackage{amsmath}
\usepackage{amsfonts}
\usepackage{latexsym}
\usepackage[all]{xy}

\CompileMatrices

\input{boxedeps.tex}
\SetEPSFDirectory{./}
\HideDisplacementBoxes

\SetRokickiEPSFSpecial  

\theoremstyle{definition}
\newtheorem{lemma}{Lemma}

\newtheorem{prop}{Proposition}
\newtheorem{corollary}{Corollary}

\DeclareMathAlphabet{\ams}{U}{msb}{m}{n}
\DeclareMathAlphabet{\goth}{U}{euf}{m}{n}

\newfont{\calli}{cmsy10}
\newcommand{\calO}{\mbox{\calli\symbol{79}}}

\def\isom{\mbox{Isom}}

\title{Constructing hyperbolic manifolds}
\author{B. Everitt \and C. Maclachlan} 
\address{Department of Mathematical Sciences, University of Durham, Durham DH1 3LE, England}
\email{brent.everitt@durham.ac.uk}
\address{Department of Mathematical Sciences, University of Aberdeen, Aberdeen AB24 3UE, Scotland}
\email{cmac@maths.abdn.ac.uk}

\subjclass{57M50}

\begin{document}


\begin{abstract}
The Coxeter simplex with symbol 
\BoxedEPSF{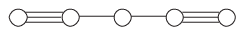 scaled 1000} 
is a compact hyperbolic
4-simplex and  the related Coxeter group $\Gamma$ is a discrete subgroup of $\text{Isom}(\ams{H}^4)$.
The Coxeter simplex with symbol 
\BoxedEPSF{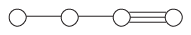 scaled 1000}
is a spherical 3-simplex, and the related 
Coxeter group $G$ is the group of symmetries of the regular 120-cell.
Using the geometry of the regular 120-cell, Davis \cite{Davis85} constructed an 
epimorphism $\Gamma\rightarrow G$ whose kernel $K$ was torsion-free, thus
obtaining a small volume compact hyperbolic 4-manifold $\ams{H}^4/K$.

In this paper 
we show how to obtain representations $\Gamma\rightarrow G$ of Coxeter groups
$\Gamma$ acting on $\ams{H}^n$ to certain classical groups $G$. We determine when the kernel $K$ of
such a homomorphism is torsion-free and thus $\ams{H}^n/K$ is a hyperbolic $n$-manifold. 
As an example, this is applied to the two groups described above,
with $G$ suitably interpreted as a classical group. Using this, further
information on the quotient manifold is obtained.
 
\end{abstract}

\maketitle


\vspace*{-2em}

\section{Introduction}

Let $M^n$ be an $n$-dimensional hyperbolic manifold, that is, an $n$-dimensional Riemannian manifold of
constant sectional curvature $-1$. Thus  $M^n$ is isometric to a quotient space $\ams{H}^n/K$
of $\ams{H}^n$ by the free action of a discrete group $K\cong\pi_1(M^n)$ of 
hyperbolic isometries.

This paper presents a method of constructing such groups $K$ as the kernels of
representations $\Gamma\rightarrow G$ of  hyperbolic Coxeter groups
$\Gamma$ into  finite classical groups. The homomorphisms arise by first
representing $\Gamma$ as a subgroup of the orthogonal group of a  quadratic space over
 a number field $k$ which preserves a lattice.  Then reducing
modulo a prime ideal in the ring of integers in the number field yields a 
representation into a finite classical group given as an orthogonal group
of a quadratic space over a finite field.

Such kernels $K$ will act freely if and only if they are torsion free. The volume of the 
resulting manifold 
$M^n=\ams{H}^n/K$ will be $N\times\text{Vol}(P)$, where $N$ is the order of the image 
in the finite classical group and $P$ is the polyhedron defining the 
Coxeter group $\Gamma$. Starting from a suitable Coxeter group $\Gamma$, the method
yields infinitely many examples of manifolds. There has been some interest lately
in constructing small volume examples when $n \geq 4$ (see \cite{Ratcliffe98a,Ratcliffe98b}).  
In dimension 4, the compact Davis manifold \cite{Davis85} is constructed
by a geometric  technique using the existence of a regular compact 120-cell in $\ams{H}^4$, 
which has volume $26 \times 4 \pi^2/3$. As an application of our method, 
we construct a compact 4-manifold $M_0$  of the same volume which turns out to be 
isometric to the Davis manifold. With the help of computational techniques,
our method gives additional information, producing a 
presentation for the fundamental group from which we obtain that $H_1(M_0) = 
\ams{Z}^{24}$.

\section{Finite representations of hyperbolic Coxeter groups}\label{method}

Consider an $(n+1)$-dimensional real space $V$
equipped with a quadratic form $q$ of signature $(n,1)$. Thus with respect to
an orthogonal basis,
\begin{equation}
q({\bf x})=-x_{n+1}^2+\sum_{i=1}^n x_i^2.
\end{equation}
The quadratic form $q$ determines a symmetric bilinear form 
\begin{equation}
B({\bf x},{\bf y}) := q({\bf x} + {\bf y}) - q({\bf x}) - q({\bf y}),
\end{equation}
with $B({\bf x},{\bf x}) = 2q({\bf x})$. The Lobachevski (or hyperboloid) model of $\ams{H}^n$ is 
the positive sheet of the sphere of unit imaginary radius in $V$ (see \cite{Ratcliffe94}).
Equivalently, we take the projection of the open cone $C = \{ {\bf x} \in V \mid q({\bf x}) < 0 ~{\rm
and~}x_{n+1} > 0 \}$ with the induced form.  

The isometries of the model are the positive $(n+1)\times (n+1)$ Lorentz matrices, that
is, the orthogonal maps of $V$, with respect to $q$, that map $C$ to itself. 

A hyperplane in $\ams{H}^n$ is the image of the intersection with $C$ of a Euclidean hyperplane in $V$. 
Each hyperplane is the projective image of the orthogonal
complement ${\bf e}^{\bot}=\{ {\bf x} \in C \mid B({\bf x},{\bf e}) = 0 \}$ of a vector ${\bf e}$
with $q({\bf e})>0$. Such a vector is said to be space-like, and it is convenient to normalise so that
$q({\bf e})=1$. The map $r_{{\bf
e}}:V\rightarrow V$ defined by
$$
r_{{\bf e}}({\bf x})={\bf x}-B({\bf x},{\bf e}){\bf e},
$$
when restricted to $\ams{H}^n$, is the reflection in the hyperplane corresponding to ${\bf e}$.

A polyhedron $P$ in $\ams{H}^n$ is the intersection of a finite collection of half-spaces, that is, the
image of
$$
\Lambda=\{ {\bf x} \in C \mid B({\bf x},{\bf e}_i) \leq 0, i = 1,2, \ldots ,m \},
$$
for some space-like vectors ${\bf e}_i$.
The intersections of the hyperplanes ${\bf e}_i^{\bot}$ with $P$ are the faces of the polyhedron.
The dihedral angle $\theta_{ij}$ subtended by two
intersecting faces of $P$ is determined by
$-2\cos\theta_{ij}=B({\bf e}_i,{\bf e}_j)$. On the other hand, non-intersecting faces of $P$ have a
common perpendicular geodesic of length $\eta_{ij}$, where
$-2\cosh\eta_{ij}=B({\bf e}_i,{\bf e}_j)$. 

All of this information is encoded in the Gram matrix $G(P)$ of $P$, an $m\times m$ matrix with
$(i,j)$-th entry $a_{ij}=B({\bf e}_i,{\bf e}_j)$. 
Let $\Gamma$ be the group generated by the reflections $r_i:=r_{{\bf e}_i}$
in the faces of $P$ so that  $\Gamma$ is a subgroup of the isometry group $\isom(\ams{H}^n)$. Moreover, $\Gamma$ is
discrete exactly when all the dihedral angles $\theta_{ij}$ of $P$ are integer submultiples
$\pi/n_{ij}$ of
$\pi$ \cite{Ratcliffe94}, and in this case, $\Gamma$ is a hyperbolic Coxeter group. The polyhedron $P$ is depicted by means of
its Coxeter symbol, with a node for each face, two nodes joined by $n-2$ edges when the corresponding faces
subtend a dihedral angle of $\pi/n$ and other pairs of nodes joined by an 
edge labelled with the geodesic length between the faces. We use the Coxeter symbol to denote
both $P$ and the group $\Gamma$ arising from it. 

Using the Lobachevski model, such a hyperbolic Coxeter group $\Gamma$ is a subgroup
of $O(V,q)$. Using this, \cite{Vinberg} gave necessary and sufficient 
conditions for such a group to be arithmetic. We adopt Vinberg's method to 
conveniently describe the groups $\Gamma$, although they need not be arithmetic.

We first give this method and some general notation which we will use throughout. Let $M$
be a finite-dimensional space over a field $F$. Equipped with a quadratic form $f$
which induces a symmetric bilinear form (as for instance in (1) and (2)), $M$ is a 
quadratic space over $F$. The group $O(M,f)$ of orthogonal maps consists of 
linear transformations $\sigma : M \rightarrow M$ such that $f(\sigma(m)) = f(m)$ 
for all $m \in M$.

Consider the Gram matrix $G(P) = [a_{ij}]$, 
and for any $\{i_1,i_2,
\ldots ,i_r\} 
\subseteq
\{1,2,
\ldots , m \}$, define
$$b_{i_1i_2 \cdots i_r} = a_{i_1i_2}a_{i_2i_3} \cdots a_{i_ri_1},$$
and let $k = \ams{Q}(\{b_{i_1i_2 \cdots i_r} \})$. 
Take the space-like vectors in $V$ defined by
$${\bf v}_{i_1i_2 \cdots i_r} = a_{1i_1}a_{i_1i_2} \cdots a_{i_{r-1}i_r}{\bf e}_{i_r}.$$
 Let $M$ be the $k$-subspace of $V$ spanned by the 
${\bf v}_{i_1i_2 \ldots i_r}$.
A simple calculation gives 
\begin{equation}\label{reflection}
r_i({\bf v}_{i_1i_2 \cdots i_r}) = {\bf v}_{i_1i_2 \cdots i_r} - {\bf v}_{i_1i_2 
\cdots i_r i},
\end{equation}
and
\begin{equation}\label{product}
B({\bf v}_{i_1i_2 \cdots i_r},{\bf v}_{j_1j_2 \cdots j_s}) = 
b_{1i_1\cdots i_rj_s\cdots j_1}.
\end{equation}
Thus, $M$ is a quadratic space over $k$ under the restriction of $q$,
and from (3) and (4)
$$B(r_i({\bf v}_{i_1i_2\cdots i_r}),r_i({\bf v}_{j_1j_2 \cdots j_s}))= B({\bf v}_{i_1i_2 \cdots
i_r},{\bf v}_{j_1j_2 \cdots j_s}).$$ 
It follows that  $\Gamma\rightarrow O(M,q)$. 

\begin{lemma}
$M$ is an $(n+1)$-dimensional space over $k = \ams{Q}(\{b_{i_1i_2 \cdots i_r} \})$.
\end{lemma}

\begin{proof}[{\bf Proof}]
If $P$ has finite volume then the vectors ${\bf e}_i$ span $V$ and the 
Gram matrix is indecomposable
\cite{Vinberg}. So for each $i$, there is a $j\not= i$ such that $a_{ij}\not= 0$. 
Successively choose indices $1=i_0,i_1,\ldots$ such that the $i_{k}$-th row contains a non-zero
entry in the $(i_{k+1})$-st column, for $k\geq 1$. 
We can ensure that the $i_k$ are distinct. For, if the only non-zero
entries of the $k$-th row are those in the columns with indices $1,i_1,\ldots,i_k$, throw away $i_k$ 
and go back to the $(i_{k-1})$-st row to rechoose a different column. Eventually, by discarding and 
moving backwards, we must be able to rechoose, in the $i_j$-th row, an index different from all the 
$i_{j+1},\ldots,i_k$ discarded. Otherwise, $\{{\bf e}_1,\ldots,{\bf e}_k\}$ are orthogonal to the other
basis vectors, contradicting indecomposability.
In this way we must arrive at a sequence $1=i_0,i_1,\ldots,i_{m-1}$ of length $m$. Hence, for any
$i$, ${\bf e}_i={\bf e}_{i_k}$ for some $i_k$, and ${\bf v}_{i_1\cdots i_k}=a_{1i_1}a_{i_1i_2}\cdots
a_{i_{k-1}i_k}{\bf e}_{i_k}$ with coefficient non-zero.
Thus, the
vectors ${\bf v}_{i_1i_2 \cdots i_r}$ span $V$ over $\ams{R}$ and hence $M$ is
$(n+1)$-dimensional over
$\ams{R}$. Now, if $\{{\bf v}_1,\ldots,{\bf v}_{n+1}\}$ is an $\ams{R}$-basis for $M$ and ${\bf v}=\sum x_i{\bf v}_i\in M$,
then the system of equations $B({\bf v},{\bf v}_j)=\sum x_i B({\bf v}_i,{\bf v}_j)$ has a unique solution,
since the matrix with $(i,j)$-th entry $B({\bf v}_i,{\bf v}_j)$ is invertible. But the solutions $x_i\in k$, since
$B({\bf u},{\bf v})\in k$ for all ${\bf u},{\bf v}\in M$. Thus $\{{\bf v}_1,\ldots,{\bf v}_{n+1}\}$
is a $k$-basis for $M$.
\end{proof}

We make a number of simplifying assumptions which hold for many examples.
Suppose that $k$ is a number field and let $\calO$ denote the ring of integers
in $k$. Suppose furthermore that all $b_{i_1 \cdots i_r} \in \calO$.
Finally let $N$ be the $\calO$-lattice in $M$ spanned by the elements
${\bf v}_{i_1 \cdots i_r}$ and assume that $N$ is a free $\calO$-lattice.
This will hold, in particular, when $\calO$ is a principal ideal domain.

By (3) above, $N$ is invariant under $\Gamma$ so that 
\begin{equation}
\Gamma \subset O(N,q) := \{ \sigma \in O(M,q) \mid \sigma(N) = N \}.
\end{equation}
With the restriction of $q$, $N$ is a quadratic module over $\calO$. If 
$\goth{P}$ is any prime ideal in $\calO$, let $\bar{k} = \calO/\goth{P}$.
Reducing modulo $\goth{P}$, we obtain a quadratic space $\bar{N}$ over $\bar{k}$
with respect to $\bar{q}$ and an induced map $\Gamma \rightarrow O(\bar{N}, \bar{q})$.

The groups $O(\bar{N},\bar{q})$ are essentially the finite classical groups
referred to earlier. However, the quadratic space $(\bar{N},\bar{q})$ may not be a 
regular quadratic space, in which case we must factor out the radical to obtain
a regular quadratic space (see Section 4 below). This will occur if the discriminant of $\bar{N}$
is zero. Since the discriminant of $\bar{N}$ is the image in $\bar{k}$ of the 
discriminant of $N$ this will only occur for finitely many prime ideals $\goth{P}$.

We now attend to the matter of when the kernel of a representation of $\Gamma$
is torsion-free. In certain  circumstances, this can be decided arithmetically
using a small variation of a result of Minkowski (see for example \cite[page 176]{Newman}).

\begin{lemma}
Let $k$ be a quadratic number field, whose ring of integers $\calO$ is a principal
ideal domain. Let $p$ be a rational prime. 
Let $\alpha \in \calO$ be such that $\alpha \not{\mid}~ 2$, and, if 3 is ramified
in the extension $k \mid \ams{Q}$, then $\alpha \not{\mid}~ 3$. If $A \in GL(n,\calO)$
is such that $A^p = I$ and $A \equiv I ({\rm mod}~\alpha)$, then $ A = I$.
\end{lemma}

\begin{proof}[{\bf Proof}]
Suppose $A \neq I$ so that $A = I + \alpha E$ where $E \in M_n(\calO)$
and we can take the g.c.d. of the entries of $E$ to be 1. From 
$(I + \alpha E)^p = I$ we have
$$ pE + \frac{p(p-1)}{2} \alpha E^2 \equiv 0\,({\rm mod}~\alpha^2).$$
Thus $p E \equiv 0\,({\rm mod}~\alpha)$ and so $p \equiv 0\,({\rm mod}~\alpha)$.
Since $\alpha \not{\mid}~ 2$, $p$ is odd. Suppose $p$ is unramified in the extension
$k \mid \ams{Q}$, then  either $p = \alpha$ 
or $p = \alpha \alpha'$ with $\alpha' \in \calO$ and $(\alpha, \alpha') = 1$.
So $pE \equiv 0\,({\rm mod}~\alpha^2)$ so $E \equiv 0\,({\rm mod}~\alpha)$.
This is a contradiction.

Now suppose that $p$ is ramified. Then $p = u \alpha^2$ where $u \in \calO^{*}$
and by assumption $ p \neq 3$. Expanding as above, but to three terms, gives
$$u\alpha^2 E + u \alpha^2 \frac{p-1}{2} \alpha E^2 + u \alpha^2 \frac{(p-1)(p-2)}{6} 
\alpha^2 E^3 \equiv 0\,({\rm mod}~\alpha^3).$$
This yields the contradiction $E \equiv 0\,({\rm mod}~\alpha)$.
\end{proof}

\begin{corollary}
If $\alpha \not{\mid}~ 2$  and if 3 is ramified in $k \mid \ams{Q}$, $\alpha \not{\mid}~ 3$, then 
the kernel of the mapping on $GL(n,\calO)$ induced by reduction $({\rm mod}~\alpha)$
is torsion-free.
\end{corollary}

More generally, 
a geometrical argument allows us to determine when {\em any\/} representation has torsion-free
kernel, albeit by expending a little more effort.
Suppose $v\in P$ is a vertex of the polyhedron
$P$, and $\Gamma_v$ is the stabiliser in $\Gamma$ of $v$. For $P$ of finite volume, $v$ is either in $\ams{H}^n$ or
on the boundary, and $v$ is called finite or ideal respectively.
We have the following ``folk-lore'' result,

\begin{lemma}\label{torsion_a}
Suppose $\Gamma$ is a discrete group generated by reflections in the faces of some polyhedron $P$ in
$n$-dimensional Euclidean space
$\ams{E}^n$ or $n$-dimensional hyperbolic space $\ams{H}^n$. If
$\gamma\in\Gamma$ is a torsion element, then for some vertex $v\in P$, $\gamma$ is
$\Gamma$-conjugate to an element of $\Gamma_v$.
\end{lemma}


Notice that in the situation described in the lemma, a Coxeter symbol for $\Gamma_v$ is obtained in
the following way: take the sub-symbol of $\Gamma$ with nodes (and their mutually incident
edges) corresponding to faces of $P$ containing $v$. For brevity's sake, when we say torsion
element from now on, we will mean {\em non-trivial\/} torsion element.

\begin{corollary}\label{}
The kernel of a representation $\Gamma\rightarrow G$ is torsion-free exactly when every torsion
element of every vertex stabiliser $\Gamma_v$ has the same order as its image in $G$.
\end{corollary}

At this point the situation bifurcates into two cases:
if $v\in P$ is a finite vertex, then $\Gamma_v$ is isomorphic to a discrete group acting on
the $(n-1)$-sphere $S^{n-1}$ centered on $v$, hence is finite. Thus, the conditions of the corollary
are satisfied exactly when $\Gamma_v$ and its image in $G$ have the same order.

If $v$ is ideal, then consider a horosphere $\Sigma$ based at $v$, and restrict the action
of $\Gamma_v$ to $\Sigma$. Then $\Sigma$ is isometric to an $(n-1)$-dimensional Euclidean space
$\ams{E}^{n-1}$, and $\Gamma_v$ acts on it discretely with fundamental region $P'$, the
intersection with $\Sigma$ of $P$. Any torsion element of $\Gamma_v$ is then $\Gamma_v$-conjugate
by Lemma \ref{torsion_a} to the stabiliser (in $\Gamma_v$!) of a vertex $v'$ of $P'$. Write
$\Gamma_{v,v'}$ for this stabiliser, and observe that it is isomorphic to a discrete group acting
on the $(n-2)$-sphere $S^{n-2}$ in $\ams{E}^{n-1}$, centered on $v'$, and hence is also finite.  
The conditions
of the corollary are satisfied exactly when for each $v'\in P'$, the group $\Gamma_{v,v'}$ and its
image in $G$ have the same order. 

Summarising,

\begin{prop}
Suppose $\Gamma$ is a hyperbolic Coxeter group generated by reflections in the faces of a
polyhedron $P$ as above. For each finite vertex $v$ of $P$, take the stabiliser $\Gamma_v$. For
each ideal vertex, take the stabilisers $\Gamma_{v,v'}$ for each vertex $v'$ of the Euclidean
polyhedron $P'$. Then $\text{kernel}(\Gamma\rightarrow G)$ is torsion-free if and only if each
such $\Gamma_{v}$ and $\Gamma_{v,v'}$ has the same order as its image in $G$.
\end{prop}

It is an elementary process to verify the conditions of the proposition. For, each
vertex stabiliser is a finite spherical reflection group of some lower dimension, hence from the
well-known list (see \cite{Humphreys90}, Section 2.11 for their orders). To find the orders of
their images in $G$, the computational algebra package {\sc Magma} is enlisted.  


\section{Polyhedra in $\ams{H}^n$}

Let $P$ be a polyhedron in $\ams{H}^n$, thus the image of
$$\Lambda=\{ {\bf x} \in C \mid B({\bf x},{\bf e}_i) \leq 0, i = 1,2, \ldots ,m \},$$
for some space-like vectors ${\bf e}_i$.
On occasion, a connected union of several copies of $P$ will yield another polyhedron of interest.
In particular, we may want to glue copies of $P$ onto its faces using some of the reflections
$r_i$ as glueing maps.

\begin{figure}
\begin{center}
\begin{tabular}{ccc}
\xymatrix{\Delta_1\,{\BoxedEPSF{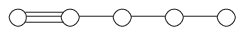 scaled 1000}}\\
\Delta_3\,{\BoxedEPSF{4.manifold.fig3.eps scaled 1000}}}&
\xymatrix{\Delta_2\,{\BoxedEPSF{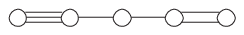 scaled 1000}}\ar@{-}[d]_-{2}\\
\Delta_4\,{\BoxedEPSF{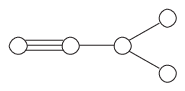 scaled 1000}}}&
$\Delta_5\,\BoxedEPSF{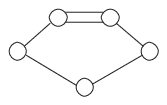 scaled 1000}$\\
\end{tabular}
\end{center}\caption{}\label{tree}
\end{figure}

\begin{lemma}
If $r_i$ is a reflection in a face of $P$, then
$$\Lambda \cup r_i(\Lambda) = \Lambda' := \{ {\bf x} \in C \mid B({\bf x},{\bf e}_j)\text{ and }
B({\bf x},r_i({\bf e}_j)) \leq 0, \text{ for all }j\not= i\}.$$
\end{lemma}

\begin{proof}[{\bf Proof}]
If ${\bf x}\in \Lambda'$, then either $B({\bf x},{\bf e}_i)\leq 0$, in which case ${\bf x}\in\Lambda$,
or $B({\bf x},{\bf e}_i)>0$, in which case $B({\bf x},r_i({\bf e}_i))\leq 0$, hence ${\bf x}\in r_i(\Lambda)$.
Conversely, if ${\bf x}\in\Lambda$, then $B({\bf x},{\bf e}_j)\leq 0$ for all $j$. If $j\not= i$, then
$$
B({\bf x},r_i({\bf e}_j))=B({\bf x},{\bf e}_j)-B({\bf x},{\bf e}_i)B({\bf e}_i,{\bf e}_j)\leq 0,
$$
since all three terms are $\leq 0$. A similar argument deals with the ${\bf x}\in r_i(\Lambda)$.
\end{proof}

We illustrate the lemma by considering the situation in four dimensions. In particular, if 
$P$ is a compact simplex it has Coxeter symbol one of the five depicted in Figure \ref{tree} (see
\cite{Humphreys90}, Section 6.9). In fact, and this explains the idiosyncratic numbering,
$\text{Vol}(\Delta_i)<\text{Vol}(\Delta_j)$ if and only if $i<j$ (see \cite{Kellerhals98}). Suppose
the nodes of
$\Delta_2$, read from left to right, correspond to hyperplanes ${\bf e}_i^{\bot}$ for
$i=1,\ldots,5$. If
$r_5=r_{{\bf e}_5}$, we have
$$
\Delta_2\cup r_5(\Delta_2)=\{{\bf x}\in C\,|\,B({\bf x},{\bf e}_i)\leq 0, i=1,\ldots,4,\text{ and
}B({\bf x},r_5({\bf e}_4))\leq 0\},
$$
since $r_5({\bf e}_i)={\bf e}_i$ for $i=1,2$ and $3$. Now, $B({\bf e}_3,r_5({\bf
e}_4))=-2\cos\pi/3$ and $B({\bf e}_4,r_5({\bf e}_4))=-2\cos\pi/2$, so $\Delta_2\cup r_5(\Delta_2)$
is a simplex with Coxeter symbol $\Delta_4$. Thus, if $\Gamma_i$ is the group generated by the
reflections in the faces of $\Delta_i$, we have that $\Gamma_4$ has index two in $\Gamma_2$. By
comparing the volumes of the simplices using the results of \cite{Kellerhals98}, the only
other possible inclusions are $\Gamma_4$ and $\Gamma_3$ as subgroups of indices 17 and 26
respectively in $\Gamma_1$. But a low index subgroups procedure in {\sc Magma} shows that
$\Gamma_1$ has no subgroups of these indices. Thus Figure \ref{tree} is a complete picture of the
possible inclusions. 


\section{An example}

In this section, we apply our method in dimension 4 starting with the Coxeter
simplex $\Delta_3$ and related group $\Gamma_3$ described above.
If $P$ is a finite volume Coxeter polyhedron in $\ams{H}^4$, then 
${\rm vol}(P) = \chi(P) 4 \pi^2/3$ where $\chi(P)$ is the Euler characteristic 
of $P$ (see \cite{Gromov82}), which coincides with the Euler characteristic of the associated
group. This is readily computed from the Coxeter symbol \cite[page 250]{Brown}, \cite{Chiswell92}.
For $\Delta_3$, the Euler characteristic is 26/14400. The vertex stabilisers are
\BoxedEPSF{4.manifold.fig6.eps scaled 1000}, $\ams{Z}_2\times \BoxedEPSF{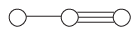 scaled
1000}$, $\BoxedEPSF{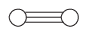 scaled 1000}\times \BoxedEPSF{4.manifold.fig8.eps scaled
1000}$, $\BoxedEPSF{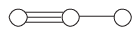 scaled 1000}\times \ams{Z}_2$ and 
\BoxedEPSF{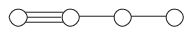 scaled 1000}, 
having orders $14400, 240,100, 240$ and $14400$ (see \cite{Humphreys90},
Section 2.11).
Thus the minimum index any torsion free subgroup of $\Gamma_3$ can have
is 14400, and we show that there is a normal torsion free subgroup of 
precisely this index. The corresponding manifold
then has Euler characteristic 26, making it the same volume as the Davis manifold \cite{Davis85}.
Indeed, it has been shown in \cite{Ratcliffe99} that $\Gamma_3$ has a unique torsion-free
normal subgroup of index 14400. It follows that this manifold is the Davis manifold.

Note that 
$$ G(\Delta_3) = \begin{pmatrix} 2 & -2 \cos \pi/5 & 0 & 0 & 0 \\
              -2 \cos \pi/5 & 2 & -1 & 0 & 0 \\
              0 & -1 & 2 & -1 & 0 \\
              0 & 0 & -1 & 2 & -2 \cos \pi/5 \\
              0 & 0 & 0 & -2 \cos \pi/5 & 2 
\end{pmatrix},
$$
so $k=\ams{Q}(\sqrt{5})$ and $M$ is a 5-dimensional space over $\ams{Q}(\sqrt{5})$. Notice that all the
$a_{ij}$, and hence the $b_{i_1i_2 \cdots i_r}$, are algebraic integers. Now, ${\bf v}_2 = -(1+\sqrt{5})/2~ {\bf e}_2$, ${\bf v}_{21} = (3+\sqrt{5})/2~ {\bf e}_1$,
${\bf v}_{23} = (1+\sqrt{5})/2~ {\bf e}_3$, ${\bf v}_{234} = -(1+\sqrt{5})/2~
{\bf e}_4$, ${\bf v}_{2345} = (3+\sqrt{5})/2~ {\bf e}_5$, and since these coefficients are all
units in $\calO$, we have $N=\calO$-span of $\{{\bf e}_1,\ldots,{\bf e}_5\}$.
The criterion of \cite{Vinberg} show that $\Gamma_3$ is arithmetic.

Let ${\goth{P}}=\sqrt{5}\calO$, so that $\bar{k}$ is the field $\ams{F}_5$ of
five elements.
By Lemma 2, the kernel of $\Gamma \rightarrow O(\bar{N},\bar{q})$ is 
torsion free. Since $(1+\sqrt{5})/2 \equiv -2\,({\rm mod}~\goth{P})$, we have that
$$\bar{q}({\bf x}) = x_1^2 + 2x_1x_2 + x_2^2 - x_2x_3 + x_3^2 - x_3x_4 + x_4^2 
+ 2x_4x_5 + x_5^2.$$
We use the same letters $\{{\bf e_1,e_2, \ldots , e_5}\}$ for the basis of 
$\bar{N}$. The images of the generating reflections of $\Gamma_3$ are then 
$5 \times 5 $ matrices with entries in $\ams{F}_5$. The computational system {\sc Magma}
then shows that the group they generate has order 14400 so that the kernel has index 14400
in $\Gamma_3$ as required.

The index 14400 is too large to allow {\sc Magma} to implement the Reidemeister-Schreier
process to obtain a presentation for $K$. However, closer examination
of the image group allows this process to be implemented by splitting into two 
steps. We will deal with this now.

The bilinear form $\bar{B}$ on $\bar{N}$ is degenerate and there is a one-dimensional
radical $\bar{N}^{\bot}$ spanned by ${\bf v}_0 = {\bf e}_1 - {\bf e}_2 + {\bf e}_4
- {\bf e}_5$. Thus $\bar{N} = W \oplus \bar{N}^{\bot}$. If ${\bf w} \in W$ and 
$\sigma \in O(\bar{N}, \bar{q})$, then $\sigma({\bf w}) = {\bf w}' + t{\bf v}_0$
where ${\bf w}' \in W$ and $t \in \ams{F}_5$. The induced mapping $\bar{\sigma}$ 
defined by $\bar{\sigma}({\bf w}) = {\bf w}'$ is easily seen to be an orthogonal map 
on $W$ and we obtain a representation $\Gamma_3 \rightarrow O(W,\bar{q})$.
We now identify $O(W,\bar{q})$ as one of the classical finite groups using the 
notation in \cite{Liebeck90}. Let ${\bf g}_i, {\bf h}_i \in W$, for $i = 1,2$ 
be defined by ${\bf g}_1 = {\bf e}_1 - {\bf e}_2$, ${\bf h}_1 = -{\bf e}_1 
+ {\bf e}_2 + {\bf e}_3$, ${\bf g}_2 = {\bf e}_1 + 2{\bf e}_5$ and ${\bf h}_2 = 
-{\bf e}_1 + 2{\bf e}_5$. Then $\bar{q}({\bf g}_i) = \bar{q}({\bf h}_i) = 0$
and $\bar{B}({\bf g}_i, {\bf h}_j) = \delta_{ij}$. Thus $O(W,\bar{q}) 
\cong O_4^{+}(5)$. There is a chain of subgroups 
$$ 1 \stackrel{2}{\subset} Z \stackrel{3600}{\subset} \Omega_4^{+}(5) \stackrel{2}{
\subset} SO_4^{+}(5) \stackrel{2}{\subset} O_4^{+}(5),$$
where $Z$ is the largest normal soluble subgroup of $\Omega_4^{+}(5)$ and 
\begin{equation}
\Omega_4^{+}(5) \cong \frac{SL(2,5) \times SL(2,5)}{\langle(-I,-I)\rangle}.
\end{equation}
The image of $\Gamma_3$ is isomorphic to a subgroup of index 2 in $O_4^{+}(5)$,
different from $SO_4^{+}(5)$ and the orientation-preserving subgroup $\Gamma_3^{+}$
maps onto $\Omega_4^{+}(5)$. The target group has a normal subgroup of index 60 with quotient isomorphic to $PSL(2,5)$ and hence so does $\Gamma_3^{+}$. Using 
{\sc Magma} we find a presentation for this subgroup $K_1$ with three generators and 
nine relations. The group $K$ is then the kernel of the induced map from $K_1$ onto $SL(2,5)$.
Again using {\sc Magma} , we obtain a presentation for $K$ on 24 generators and several
pages of relations. The abelianisation of $K$ is $\ams{Z}^{24}$. This agrees with the homology 
calculations in \cite{Ratcliffe99}.

This calculation is readily carried out once the images of the generators of 
$\Gamma_3^{+}$ are identified with pairs of matrices.
We sketch the method of obtaining this description.

Let $V$ be a two dimensional space over $\ams{F}_5$ with symplectic form $f$ defined with respect to a basis 
${\bf n}_1, {\bf n}_2$ by 
$$f(\sum x_i {\bf n}_i, \sum y_i {\bf n}_i) = x_1y_2 - x_2 y_1.$$
Let $U = V \otimes V$ and define $g$ on $U$ by
$$g({\bf v}_1 \otimes {\bf v}_2, {\bf w}_1 \otimes {\bf w}_2) = f({\bf v}_1, {\bf w}_1) f({\bf v}_2,{\bf w}_2).$$
Then $g$ is a symmetric bilinear form on $U$ and $O(U,g) \cong O_4^{+}(5)$.
Note that $SL(2,5) \times SL(2,5)$ acts on $U$ by
$$ (\sigma, \tau)({\bf v} \otimes {\bf w}) = \sigma({\bf v}) \otimes \tau({\bf w}),$$
and this action preserves $g$ with $(-I,-I)$ acting trivially. This describes
the group $\Omega_4^{+}(5)$. Additionally, the mapping $\rho : U \rightarrow U$
given by $\rho({\bf v} \otimes {\bf w}) = {\bf w} \otimes {\bf v}$ also lies in
$O(U,g)$ and has determinant -1. Let $H$ be the subgroup generated by $\Omega_4^{+}(5)$ and $\rho$. 

We identify the image of $\Gamma_3$ with $H$, by first identifying $U$ and $W$ by the linear isometry induced by 
$${\bf n}_1 \otimes {\bf n}_1 \mapsto {\bf g}_1, {\bf n}_1 \otimes {\bf n}_2 \mapsto 
{\bf g}_2, {\bf n}_2 \otimes {\bf n}_1 \mapsto - {\bf h}_2, {\bf n}_2 \otimes 
{\bf n}_2 \mapsto {\bf h}_1.$$
It is now easy to check that the image of $r_5$ is $\rho$. Recall that 
$\Gamma_3^{+}$ is generated by $x = r_5 r_4, y = r_5 r_3, z = r_5 r_2, w = r_5 r_1$.
Now determine the images of $x,y,z,w$ as pairs of matrices in $SL(2,5) \times SL(2,5)$.

\small{}

\end{document}